\documentclass[12pt]{amsart}

\usepackage{amssymb}
\usepackage{graphicx}  
\DeclareGraphicsExtensions{.pstex,.eps}

\newcommand{\RR}{{\mathbb R}}

\theoremstyle{plain}

\newtheorem{thm}{Theorem}
\newtheorem{prop}{Proposition}[section]

\newtheorem{lem}[prop]{Lemma}

\theoremstyle{definition}

\numberwithin{equation}{section}

\def\squarebox#1{\hbox to #1{\hfill\vbox to #1{\vfill}}}

\newcommand{\p}{\partial}
\usepackage{amsxtra}

\title
[Integrability along the flow]
{Wave packet parametrices for evolutions governed by 
pdo's with rough symbols}

\author[J. Marzuola]
{Jeremy Marzuola}

\author[J. Metcalfe]
{Jason Metcalfe}

\author[D. Tataru]
{Daniel Tataru}

\address{Mathematics Department, University of California \\
Evans Hall, Berkeley, CA 94720-3840, USA}

\begin{document}    
   
\maketitle

\section{Introduction}

In this article we consider evolution equations of the form
\begin{equation} 
\left\{ \begin{array}{c}
(D_t + a^w(t,x,D)+ib^w(t,x,D)) u  =  f, \ \text{ in } \RR^+ \times \RR^n \\ \\
u(0)  = u_0, \  \text{ in } \RR^n
\end{array} \right.
\label{cp} 
\end{equation}
where $a(t,x,\xi)$ and $b(t,x,\xi)$  are real symbols which are continuous in $t$ and
smooth with respect to $x$ and $\xi$. 

The operator $a^w(t,x,D)$ is selfadjoint; if $b=0$ then this
formally guarantees that the above evolution is $L^2$ well-posed and
the corresponding evolution operators $S(t,s)$ are $L^2$ isometries.
The $b^w$ term rougly contributes to the growth or decay of energy
along the flow, depending on whether $b$ is negative or positive.

We are interested in the phase space localization 
properties of the evolution operators $S(t,s)$. These are best
described in terms of the Bargman transform,
\begin{equation}
(T f)(x,\xi) = 
2^{-\frac{n}{2}} \pi^{-\frac{3n}{4}} \int e^{-\frac{1}{2}(x-y)^2} e^{i\xi(x-y)}
f(y) \ dy, 
\end{equation}
which is an isometry from $L^2(\RR^n)$ to the subspace of
$L^2(\RR^{2n})$ of functions satisfying the Cauchy-Riemann type relation
\begin{equation}
i \p_\xi Tf = (\p_x -i\xi) Tf.
\label{cr} \end{equation}
The inversion formula is
\begin{equation}
f(y)= 2^{-\frac{n}{2}} \pi^{-\frac{3n}{4}} \int 
e^{-\frac{1}{2}(x-y)^2}  e^{i\xi(y-x)} (T f)(x,\xi) \ dx d\xi\,.
\label{inversion}\end{equation}

Then one would like to describe the phase space localization of
$S(t,s)$ relative to the Hamilton flow corresponding to
\eqref{cp}. This is given by
\begin{equation}
\left\{ \begin{array}{c} \dot x = a_\xi(t,x,\xi) \cr \dot \xi =
    -a_x(t,x,\xi) \end{array} \right.
\label{hf}\end{equation}
We denote by $\chi(t,s)$ the corresponding family of cannonical
transformations, and by 
\[
t \to (x^t,\xi^t)
\]
the trajectories of the Hamilton flow.

This problem has already been considered in \cite{T},
\cite{KT}. There, they the class
$S^{0,(k)}_{00}$ of symbols which satisfy the bounds
\begin{equation}
|\p_x^\alpha \p_\xi^\beta a(x,\xi)| \leq c_{\alpha\beta}, \qquad |\alpha|+|\beta|
\geq k.
\label{(k)}\end{equation}
The main result has the form

\begin{thm}\cite{T},\cite{KT}
Assume that the symbol $a(t,x,\xi)$ satisfies $a(t,x,\xi) \in
S^{0,(2)}_{00}$ uniformly with respect to $t$. Then

a) The Hamilton flow is bilipschitz.

b) The kernel $K(t,s)$ of the phase space operator $T^* S(t,s) T$
decays rapidly away from the graph of the Hamilton flow,
\begin{equation}
| K(t,x,\xi,s,y,\eta)| \lesssim (1+ |(x,\xi) - \chi(t,s)(y,\eta)|)^{-N}.
\end{equation}
\end{thm}

However, for applications to nonlinear evolution equations
one would like to relax the above class of symbols and 
replace uniform bounds by an integrability condition.
For instance, in the context of the wave equation related results have 
been obtained in \cite{lp} under assumptions which correspond to
replacing the $L^\infty$ bounds in \eqref{(k)} with $L^1_t
L^\infty_x$.

In this article we go one step further and restrict the time
integrability to the bicharacteristic rays. This is a much more
natural condition from the point of view of applications. One
motivation for this already appears in early works of
Mizohata~\cite{Miz1,Miz2} which is concerned with the $b^w$
type terms.  They consider the equation
\begin{eqnarray*}
Lu := \partial_t - i\Delta + \sum_{j=1}^n b_j(x) \partial_{x_j} + c(x,t) u = f(x,t), 
\end{eqnarray*}
and show that a necessary condition for $L$ to be well-posed in
$H^\infty$ is the bound
\begin{equation} 
\label{co1}
\sup_{x \in \RR^{n}, \omega \in \mathbb{S}^{n-1}, R > 0} \left|Im
  \int_{0}^{R} b_{1}(x+r \omega) \cdot \omega dr\right | < \infty. 
\end{equation}
On the other hand, a slightly stronger version of \eqref{co1} was
shown to be sufficient for $L^2$ wellposedness in \cite{BT}.

In the case where $\Delta$ is replaced by the variable coefficient operator $a_{jk} (x,t) \partial_j \partial_k$, then a natural extension of \eqref{co1} is
\begin{eqnarray*}
\sup_{x \in \RR^{n}, |\xi| = 1, R > 0} \left|Im
  \int_{0}^{R} b_{1}(x^t(t,x,\xi)) \cdot \xi^t (t,x,\xi) dr\right | < \infty.
\end{eqnarray*}

Another motivation for this work comes from the study of 
general quasilinear Schr\"odinger equations.  In \cite{KPV} and \cite{KPRV}, well-posedness is established in highly regular Sobolev spaces
by using estimates for the corresponding linear equation.

 Given a symplectic flow $\chi$ in $\RR \times \RR^{2n}$ we introduce the 
symbol class $S^{(k)}{L^1_\chi}$ of symbols $q$, which are smooth in
$(x,\xi)$, continuous in $t$ and satisfy
\begin{equation}
  \sup_{x,\xi}  \int_{0}^1 |\p_x^\alpha \p_\xi^\beta q(t,\chi(t,0)(x,\xi))|
  dt \leq c_{\alpha\beta}, \qquad |\alpha|+|\beta| \geq k.
\label{(kl1)}\end{equation}

Then our condition for the symbol $a$ is implicit, namely $a \in
S^{(2)}{L^1_\chi}$ where $\chi$ is the Hamilton flow of $a$ defined by
\eqref{hf}. For the symbol $b$ we will assume that $b \in
S^{(1)}{L^1_\chi}$. Given such $a$ and $b$ we introduce the notation
\begin{equation}
  \kappa_N = \max_{ 2 \leq |\alpha|+|\beta| \leq N} c_{\alpha\beta}^a
  + \max_{ 1 \leq |\alpha|+|\beta| \leq N} c_{\alpha\beta}^b, \qquad
  \kappa_0 = \max_{|\alpha|+|\beta| = 2} c_{\alpha\beta}^a
  + \max_{|\alpha|+|\beta| =1} c_{\alpha\beta}^b
\label{cm}\end{equation}
where $c_{\alpha\beta}^a$ and $c_{\alpha\beta}^b$ are as in
\eqref{(kl1)} corresponding to the symbols $a$ and $b$.

The other important parameter in our analysis corresponds 
to \eqref{co1}. We set
\begin{equation}
M = \sup_{x,\xi} \sup_{0 \leq t_0 \leq t_1}  \int_{t_0}^{t_1} b(t,x^t,\xi^t) dt
\label{M}\end{equation}
and assume that $M$ is finite. Then our main result is

\begin{thm}
  a) Assume that the symbol ${a}(t,x,\xi)$ satisfies $a(t,x,\xi) \in
  S^{(2)}L^1_\chi$. Then the Hamilton flow defined by \eqref{hf} is
  globally well defined and bilipschitz.

  b) Assume in addition that $b$ is a  symbol in $S^{(1)}{L^1_\chi}$
so that $M$ given by \eqref{M} is finite and the following relation
holds for some large $N$:
\[
e^{2M} \kappa_0 \kappa_{4N} \ll 1.
\]
  Then the kernel $K(t,s)$ of the phase space operator $T^* S(t,s) T$
  decays rapidly away from the graph of the Hamilton flow,
\begin{equation}
| K(t,x,\xi,s,y,\eta)| \lesssim (1+ |(x,\xi) - \chi(t,s)(y,\eta)|)^{-N}.
\label{K}\end{equation}
\label{maint}\end{thm}

We remark that the smallness in part (b) can be replaced by an
equiintegrability condition for the second order derivatives of $a$ 
along the flow,
\begin{equation}
  \lim_{h \to 0} \sup_{x,\xi,t_0}  \int_{t_0}^{t_0+h} |\p^2 a(x^t,\xi^t)|
  dt =0.
\label{(kl2)}\end{equation}
Then one can take $N$ arbitrary, apply the theorem on sufficiently
small time intervals and iterate.

\medskip
\noindent
{\sc Acknowledgments.} 
The work of the second author was supported in part by an NSF
postdoctoral fellowship, and that of the first and third authors by
NSF grants DMS0354539 and DMS0301122.

\section{The Hamilton flow}

In this section we prove that if $a \in S^{(2)}L^1_\chi$ then the
Hamilton flow for \eqref{cp} is well defined and bilipschitz. We first
prove the bilipschitz property locally, and then use it to show that
the flow is globally well defined in the entire time interval.  

Thus begin with $(x^0,\xi^0) \in \RR^{2n}$. Then there exists some
time $t_0 > 0$ and a ball $B$ centered at $(x^0,\xi^0)$ so that we can
solve \eqref{hf} with initial data $(y^0,\eta^0) \in B$. Since, by standard
ODE results the flow maps $\chi(t,s)$ are smooth, we only need to
obtain the appropriate bounds.

For simplicity, set
\begin{eqnarray*}
\vec{p}(t,x,\xi) = \left( \begin{array}{c}
x^t (x,\xi) \\
\xi^t (x,\xi)
\end{array} \right).
\end{eqnarray*}
We have the following systems of equations:
\begin{eqnarray*}
\p_x \vec{p}(t) & = & \left( \begin{array}{c}
1 \\ 
0
\end{array} \right) + \int_0^t \left( \begin{array}{cc}
a_{\xi x}(s) & a_{\xi \xi}(s) \\
-a_{x x}(s) & -a_{\xi x}(s)
\end{array} \right) \p_x \vec{p}(s) ds \\
& = & \left( \begin{array}{c}
1 \\ 
0
\end{array} \right) + \int_0^t A(s) \p_x \vec{p}(s) ds,
\end{eqnarray*}
and a similar expression for $\p_\xi \vec{p}(t)$ for all $0 \leq t \leq 1$.

Now, let us see that the flow is Lipschitz in $x$.  Taking absolute values and applying Gronwall's inequality, we have
\begin{eqnarray*}
\|\partial_x \vec{p}(t) \| \leq e^{\int_0^t \|A(s)\| ds}.
\end{eqnarray*}
However,
\begin{eqnarray*}
\| A \| \lesssim | \partial^2 a(s)|,
\end{eqnarray*}
thus, since $a \in S^{(2)}{L^1_\chi}$, 
\begin{eqnarray*}
\|\partial_x \vec{p}(t) \| \lesssim 1,
\end{eqnarray*}
and similarly
\begin{eqnarray*}
\|\partial_\xi \vec{p}(t) \| \lesssim 1.
\end{eqnarray*}

Note that bounds on higher derivatives in $x$ and $\xi$ will follow
from the same argument.

We have shown that the derivatives of $x^t$ and $\xi^t$ are uniformly
bounded; therefore it is clear that the Hamilton flow
map $\chi: (x,\xi) \to (x^t,\xi^t)$ is uniformly bilipschitz in phase
space.

It remains to show that the Hamilton flow is globally defined.
Consider a bicharacteristic $t \to (x^t,\xi^t)$ starting at
$(x^0,\xi^0)=(y,\eta)$. This can be continued in time for as long as $
(x^t,\xi^t)$ remains finite. To prove that this is indeed the case we
 consider a one parameter family of
bicharacteristics
\[
s \to (x^t(s),\xi^t(s)), \qquad (x^s(s),\xi^s(s))=(y,\eta), \ \ t \geq s 
\]
and show that $ (x^t(s),\xi^t(s))$ is of class $C^1$ with respect
to $s$ with a uniform Lipschitz bound for $s,t$ in a compact time
interval.

Since the Hamilton flow is smooth with a bounded differential
it suffices to estimate the derivative at $s=t$. But this is given by
\[
\frac{d}{ds}  (x^t(s),\xi^t(s))_{s=t} = (a_\xi(t,y,\eta),-a_x(t,y,\eta))
\] 
and due to the continuity in $t$ it is bounded on any compact time
interval.

\section{Kernel Bounds}

We begin with the equation \eqref{cp} for $u$, and, following
\cite{T}, we derive an equation for its phase space transform $Tu$.
If the symbol $a$ is linear in $x$ and $\xi$, then we have the
straightforward identity
\[
T a^w(x,D) u = (a(x,\xi) + i (a_x \partial_\xi - a_\xi (\partial_x
-i\xi) ) )Tu.
\]
For a general symbol $a$ we denote by $a_{x,\xi}$ its linearization
around $(x,\xi)$. Then
\[
(T a^w(x,D) u)(x,\xi) = (a(x,\xi) + i (a_x \partial_\xi - a_\xi
(\partial_x -i\xi) ) ) Tu(x,\xi) + T (a -a_{x,\xi})^w(x,D) u.
\]
This provides us with an evolution equation for $v = Tu$,
\begin{equation}
  (D_t +  a + i (a_x \partial_\xi - a_\xi \partial_x ) - \xi a_\xi  +ib) v
  = E v, \qquad v(0) = v_0 = Tu_0,
\label{veq}\end{equation}
where the remainder $E$ is given by
\[
Ev(x,\xi) = T [ ( a -a_{x,\xi})^w(x,D)+ i(b-b_{x,\xi})^w(x,D)] T^* v (x,\xi),
\]
where $b_{x,\xi}$ is the zeroth order term in a Taylor expansion of $b$ about $(x,\xi)$. 
The kernel $K(t,x,\xi,0,y,\eta)$ of $TS(t,0) T^*$ is given by the 
solution $v(t,x,\xi)$ to \eqref{veq} with initial data 
\[
v_0(x,\xi) = C_n \int_{\RR^n} e^{i \xi(x-z)} e^{-\frac{(x-z)^2}2}
e^{-i \eta(y-z)} e^{-\frac{(y-z)^2}2} dz = C_n e^{-\frac{(x-y)^2}4}
e^{-\frac{(\xi-\eta)^2}4}e^{i\frac{(x-y)(\xi+\eta)}2}
\]
which decays rapidly away from $(y,\eta)$.  Then, the bound \eqref{K}
follows if we show that the equation \eqref{veq} is wellposed in
weighted $L^\infty$ spaces:

\begin{prop}
Let $a$, $b$ be as in part $b$ of Theorem~\ref{maint}.  Then for each initial data 
$v_0$ satisfying
\[
(1+ |x-y|+|\xi-\eta|)^N v_0(x,\xi) \in L^\infty
\]
there exists an unique solution $v$ to \eqref{veq} satisfying 
\[
(1+ |x-y^t|+|\xi-\eta^t|)^N v(t,x,\xi) \in L^\infty
\]
\end{prop}

\begin{proof}
The inhomogeneous equation
\begin{equation}
  (D_t +  a + i (a_x \partial_\xi - a_\xi \partial_x ) - \xi a_\xi+ib) v
  = f, \qquad v(0) = v_0 
\label{veq1}\end{equation}
is an ODE along the Hamilton flow, whose solution satisfies the bound
\[
\begin{split}
|v(t,x^t,\xi^t)| \leq& |v_0(x,\xi)|  e^{\int_{0}^t
b(\sigma,x^\sigma,\xi^\sigma) d\sigma} + \int_0^t |f(s,x^s,\xi^s)|
e^{\int_{s}^t b(\sigma,x^\sigma,\xi^\sigma) d\sigma} ds
\\  \leq& e^M \left(|v_0(x,\xi)| + \int_0^t |f(s,x^s,\xi^s)|
 ds \right).
\end{split}
\]
By the bilipschitz property of the flow this implies that
\[
\begin{split}
\|(1+ |x-y^t|+|\xi-\eta^t|)^N v(t,x,\xi)\|_{L^\infty} & \lesssim  \
e^M \Bigl(
\|(1+ |x-y|+|\xi-\eta|)^N v_0(x,\xi)\|_{L^\infty}
\\ &+ \sup_{x,\xi}  (1+ |x-y|+|\xi-\eta|)^N \int_0^1 |f( s,x^s,\xi^s)| ds \Bigr).
\end{split}
\]
Hence in order to be able to treat the right hand side term $Ev$ perturbatively and
prove the proposition it suffices to show that $E$ satisfies the estimate
\begin{equation}
\sup_{x,\xi} \Bigl[ (1+ |x-y|+|\xi-\eta|)^N \int_0^1 |Ev(
s,x^s,\xi^s)| ds \Bigr] \ll e^{-M} \|(1+ |x-y^t|+|\xi-\eta^t|)^N v(t,x,\xi)\|_{L^\infty}.
\label{ebd}\end{equation}

For this we have to understand in more detail the kernel of  $E$. It
suffices to consider the contribution of $a$ in $E$, the estimates for
the contribution of $b$ are similar but simpler.
We begin with a
simple computation for the kernel $K_q$  of $T q^w(x,D) T^*$,
\[
\begin{split}
K_q(x,\xi,x_1,\xi_1) = & 
C_n \int  e^{i \xi(x-y)} e^{-\frac{(x-y)^2}2}  e^{i\eta(y-y_1)}
q(\frac{y+y_1}2,\eta)  e^{-i \xi_1(x_1-y_1)} e^{-\frac{(x_1-y_1)^2}2} dy dy_1 d\eta.
\end{split}
\]
We change variables to
\[
z = \frac{y+y_1}2, \qquad w = \frac{y-y_1}2
\]
and integrate with respect to $w$ to obtain
\[
\begin{split}
K_q(x,\xi,x_1,\xi_1) = & 
C_n \int  e^{i \xi(x-z-w)} e^{-\frac{(x-z-w)^2}2}  e^{2i\eta w}
q(z,\eta)  e^{-i \xi_1(x_1-z+w)} e^{-\frac{(x_1-z+w)^2}2} dz dw d\eta
\\ = &C_n e^{\frac{i(x+x_1)(\xi-\xi_1)}2}
\int  e^{-\frac{(\xi+\xi_1-2\eta)^2}4}  e^{-\frac{(x+x_1-2z)^2}4}  
q(z,\eta)  e^{i \eta(x-x_1)} e^{-iz(\xi-\xi_1)} dz  d\eta.
\end{split}
\]
Integrating by parts $2N$ times with respect to $\eta$ and $z$ leads 
to
\[
|K_q(x,\xi,x_1,\xi_1)| \lesssim (1+(x-x_1)^2 +
(\xi-\xi_1)^2)^{-N} \int \! e^{-\frac{(\xi+\xi_1-2\eta)^2}8}  e^{-\frac{(x+x_1-2z)^2}8}  
\sup_{|\alpha| \leq 2N} |\partial^\alpha q(z,\eta)| dz d\eta.
\]

To estimate this we use the following 
\begin{lem}
Let $q$ be a smooth function in $\RR^n$ with $q(0) = 0$, $\nabla q(0)=
0$. Then
\begin{equation}
\int_{|x| \lesssim R} |q(x)| dx \lesssim  R^{n+1} \int_{|x| \lesssim R}
|x|^{1-n} |\partial^2 q(x)| dx 
\end{equation}
\end{lem}
\begin{proof}
  The argument follows directly from the fundamental theorem of
  calculus.  We have
\begin{eqnarray*}
  \int_{|x| \lesssim R} |q(x)| dx & \lesssim &  \int_0^R \int_{S^{n-1}} \int_0^r \int_0^{r_1} | \partial^2_r q(r_2,\omega)|  dr_2 dr_1  r^{n-1} dr d \omega \\
  & \lesssim & \int_0^R \int_{S^{n-1}} \int_0^{R} | \partial^2_r q(r_2,\omega)| r_2^{1-n} r_2^{n-1}  dr_2 r^{n} dr d \omega \\
  & \lesssim & R^{n+1} \int_{|x| \lesssim R} |x|^{1-n} |\partial^2 q| dx.
\end{eqnarray*}
\end{proof}

We apply the lemma to estimate $K = K_q$ with $q = a-a_{x,\xi} $ and $N=0$.
This gives the crude bound
\begin{equation}
|K(x,\xi,x_1,\xi_1)| \lesssim (1+|\xi-\xi_1|+|x-x_1|)^{2n+3}
\int k(x-z,\xi-\eta)
|\partial^2 a(z,\eta)| dz d\eta,
\end{equation}  
where $k$ is the integrable kernel
\[
k(x,\xi) = (|x|+|\xi|)^{1-2n} (1+|x|+|\xi|)^{-2}.
\]
Applying the same bound with large $N$ yields
\begin{equation}
|K(x,\xi,x_1,\xi_1)| \lesssim (1+|\xi-\xi_1|+|x-x_1|)^{2n+3-4N}
\int k(x-z,\xi-\eta)
\sup_{2 \leq |\alpha| \leq 4N} |\partial^\alpha a(z,\eta)| dz d\eta.
\end{equation}

Now we can get bounds on $E$:
\[
\begin{split}
 & \int_{0}^1 |E v(t,x^t,\xi^t)|dt  \lesssim
\sup_{t\in [0,1]}  \int (1+|x^t-x_1|+|\xi^t-\xi_1|)^{2n+3-2N} |v(t,x_1,\xi_1)| dx_1 d\xi_1 \\
 & \int_0^1 \left( \int k(x^t-z,\xi^t-\eta) |\partial^2 a(z,\eta)| dz d\eta
  \right)^\frac12 
 \! \left( \int k(x^t-z,\xi^t-\eta)\sup_{2 \leq |\alpha|
      \leq 4N} |\partial^\alpha a(z,\eta)| dz d\eta \right)^\frac12\!\! dt
\\ &  \hspace{1.3in} \lesssim
\sup_{t\in [0,1]}  \int (1+|x^t-x|+|\xi^t-\xi|)^{2n+3-N} |v(t,x,\xi)| dx d\xi \\
 & \left( \int k(x^t-z,\xi^t-\eta) |\partial^2 a(z,\eta)| dz d\eta dt
  \right)^\frac12  \left( \int \!k(x^t-z,\xi^t-\eta)\sup_{2 \leq |\alpha|
      \leq 4N} |\partial^\alpha a(z,\eta)| dz d\eta dt \right)^\frac12.
\end{split}
\]
In the last two integrals we change variables $(z,\eta) \to
(y^t,\zeta^t)$, which uses the boundedness proved in Section $2$. Since $\chi(t,s)$ are bilipschitz we can replace
$ k(x^t-y^t,\xi^t-\zeta^t)$ by $ k(x^0-y^0,\xi^0-\zeta^0)$.
Then we use \eqref{(kl1)} for the time integral and the 
integrability of $k$ for the $(y^0,\zeta^0)$ integral. We obtain
\begin{equation}
\begin{split}
 & \int_{0}^1 |E v(t,x^t,\xi^t)|dt  
 \lesssim \sqrt{\kappa_0 \kappa_{4N}}
\sup_{t\in [0,1]}  \int (1+|x^t-x|+|\xi^t-\xi|)^{2n+3-2N} |v(t,x,\xi)| dx d\xi,
\end{split}
\label{eint}\end{equation}
where $\kappa_0$ and $\kappa_{4N}$ are as in \eqref{cm}.

This allows us to iteratively solve the equation \eqref{veq} for
initial data $v_0 \in L^\infty$. Indeed, the above estimate implies that
\[
\sup_{x,\xi} \int_{0}^1 |E v(t,x^t,\xi^t)|dt  \lesssim \sqrt{\kappa_0 \kappa_{4N}}
\|v\|_{L^\infty}.
\]

It also allows us to solve the equation \eqref{veq} for
initial data $v_0$  in weighted $L^\infty$ spaces. Precisely, for $2N > 4n+m+3 $
we have
\begin{eqnarray*}
\sup_{x,\xi} \Bigl[ (1 + |x-x_1|+& &\negthickspace \negthickspace \negthickspace \negthickspace \negthickspace \negthickspace \negthickspace \negthickspace |\xi-\xi_1|)^m \int_0^1 |Ev(t,x^t,\xi^t)|dt \Bigr] \\
& \lesssim & 
\sup_{x,\xi} \Bigl[ (1+|x^t-x_1^t|+|\xi^t-\xi_1^t|)^m \sqrt{\kappa_0 \kappa_{4N}}
 \\
& \ & \times \sup_{t\in [0,1]}  \int (1+|x^t-y|+|\xi^t-\eta|)^{2n+3-2N} |v(t,y,\eta)| dy d\eta \Bigr] \\
& \lesssim &  \sqrt{\kappa_0 \kappa_{4N}} \| (1+|x-x_1^t|+|\xi-\xi_1^t|)^m v(t,x,\xi) \|_{L^\infty_{x,\xi}}.
\end{eqnarray*}
Due to the hypothesis in part (b) of the theorem this implies
\eqref{ebd} and concludes the proof of the proposition.
\end{proof}

\end{document}